\documentclass[fleqn]{mat01}
\usepackage{times,mathtimy,amssymb}
\begin{document}

\setcounter{page}{183}
\firstpage{183}

\newtheorem{definit}{\rm DEFINITION}

\newtheorem{theo}{Theorem}
\renewcommand{\thetheo}{\arabic{theo}}
\newtheorem{theor}[theo]{\bf Theorem}

\newtheorem{corolla}{\rm COROLLARY}

\newtheorem{rema}{Remark}

\title{Fixed point of multivalued mapping in uniform spaces}

\markboth{Duran T\"{u}rkoglu and Brian Fisher}{Multivalued mapping in uniform spaces}

\author{DURAN T\"{U}RKOGLU and BRIAN FISHER$^{*}$}

\address{Department of Mathematics, Faculty of Science and Arts,
University of Kirikkale, 71450-Yahsihan Kirikkale, Turkey\\
\noindent $^{*}$Department of Mathematical Sciences,
Faculty of Computer Sciences and Engineering,
De Montfort University, The Gateway, Leicester, LE1 9BH, England\\
\noindent E-mail: duran$_{-}$t@hotmail.com; fbr@le.ac.uk}

\volume{113}

\mon{May}

\parts{2}

\Date{MS received 31 December 2001; revised 10 July 2002}

\begin{abstract}
In this paper we prove some new fixed point theorems for multivalued
mappings on orbitally complete uniform spaces.
\end{abstract}


\maketitle

\section{Introduction}

Let $(X,\mathcal{U})$ be a uniform space. A family $\left\{ d_{i}: 
i\in I\right\}$ of pseudometrics on $X$ with indexing set $I,$ is called an
associated family for the uniformity $\mathcal{U}$ if the family 
\begin{equation*}
\beta = \left\{ V(i, \varepsilon): i\in I; \varepsilon > 0 \right\},
\end{equation*}
where
\begin{equation*}
V(i,\varepsilon) = \left\{ (x,y) : x,y 
\in X,\ \ d_{i}(x,y) < \varepsilon \right\}
\end{equation*}
is a sub-base for the uniformity $\mathcal{U}$. We may assume that $\beta$ 
itself is a base by adjoining finite intersection of members of $\beta$, if
necessary. The corresponding family of pseudometrics is called an associated
family for $\mathcal{U}$. An associated family for $\mathcal{U}$ will be
denoted by $p^{*}$. For details the reader is referred to \cite{1,3,4,5,6,7,8}.

Let $A$ be a nonempty subset of a uniform space $X$. Define
\begin{equation*}
\Delta^{*}(A) = \sup \left\{ d_{i} (x,y): x,y \in A,\,i\in I\right\},
\end{equation*}
where 
\begin{equation*}
\left\{ d_{i}:i\in I\right\} = p^{*}.
\end{equation*}
Then $\Delta^{*}$ is called an augmented diameter of $A$. Further, $A$ is
said to be $p^{*}$-{\it bounded} if $\Delta ^{*}(A) < \infty$. Let 
\begin{equation*}
2^{X} = \left\{ A: A\ \hbox{ is a nonempty, closed and } p^{*}\ \hbox{-
bounded subset of } X \right\}.
\end{equation*}
For any nonempty subsets $A$ and $B$ of $X$, define 
\begin{equation*}
d_{i}(x,A) = \inf\ \{ d_{i}(x,\,a):a\in A\}\text{, }i\in I
\end{equation*}
\begin{align*}
H_{i}(A,\,B) &= \max \left\{ \underset{a\in A}{\sup }\ d_{i}(a,\,B),\,
\underset{b\in B}{\sup }\ d_{i}(A,\,b)\right\} \\[.2pc]
&= \underset{x\in X}{\sup }\left\{ \left| d_{i}(x,\,A)-d_{i}(x,\,B)\right|
\right\}.
\end{align*}

It is well-known that on $2^{X}$, $H_{i}$ is a pseudometric, called the
Hausdorff pseudometric induced by $d_{i},\,i\in I$.

Let $(X,\mathcal{U})$ be a uniform space with an augmented associated family 
$p^{*}$. $p^{*}$ also induces a uniformity $\mathcal{U}^{*}$ on $2^{X}$
defined by the base
\begin{equation*}
\beta^{*} = \left\{ V^{*}(i,\,\varepsilon): i\in I,\,\varepsilon
>0\right\},
\end{equation*}
where 
\begin{equation*}
V^{*}(i,\,\varepsilon) = \left\{ (A,\,B):A,\,B\in 2^{X},\,H_{i} (A,\,B) <
\varepsilon \right\}.
\end{equation*}
The space $(2^{X},\,\mathcal{U}^{*})$ is a uniform space called the
hyperspace of $(X,\,\mathcal{U})$.\vspace{.4pc}

\begin{definit}$\left.\right.$\vspace{.5pc}

\noindent {\rm The collection of all filters on a given set $X$ is
denoted by $\Phi (X)$. An order relation is defined on $\Phi (X)$ by the
rule $\mathcal{F}_{1}<\mathcal{F}_{2}$ iff $\mathcal{F}_{1}\supset
\mathcal{F}_{2}$. If $\mathcal{F}^{*}<\mathcal{F}$, then
$\mathcal{F}^{*}$ is called a subfilter of $\mathcal{F}$.}\vspace{.4pc}
\end{definit}

\begin{definit}$\left.\right.$\vspace{.5pc}

\noindent {\rm Let $(X,\,\mathcal{U})$ be a uniform space defined by 
$\{d_{i}:i\in I\} = p^{*}$. If $F:X\rightarrow 2^{X}$ is a multivalued
mapping, then

\begin{enumerate}
\renewcommand{\labelenumi}{(\roman{enumi})}
\leftskip .5pc
\item $x\in X$ is called a fixed point of $F$ if $x\in Fx$;

\item An orbit of $F$ at a point $x_{0}\in X$ is a sequence
$\{x_{n}\}$ given by
\begin{equation*}
O(F, x_{0}) = \{x_{n}: x_{n}\in Fx_{n-1}, n = 1, 2, ...\};
\end{equation*}
\item A uniform space $X$ is called $F$-{\it orbitally complete} if every
Cauchy filter which is a subfilter of an orbit of $F$ at each $x\in X$
converges to a point of $X$.
\end{enumerate}}
\end{definit}

\begin{definit}$\left.\right.$\vspace{.5pc}

\noindent {\rm Let $(X,\,\mathcal{U})$ be a uniform space and let $F:X\rightarrow X$ be
a mapping. A single-valued mapping $F$ is {\it orbitally continuous} if
$\lim \,(T^{n_{i}}x)=u$ implies $\lim \,T(T^{n_{i}}x) = Tu$ for each
$x\in X$.}
\end{definit}

\section{Main results}

\begin{theor}[\!]
Let $(X,\,\mathcal{U})$ be an $F$-orbitally complete Hausdorff uniform space
defined by $\left\{ d_{i}:\,i\in I\right\} =p^{*}$ and $(2^{X},\,\mathcal{U}%
^{*})$ a hyperspace and let $F:\,X\rightarrow \,2^{X}$ be a continuous
mapping with $Fx$ compact for each $x$ in $X$. Assume that 
\setcounter{equation}{0}
\begin{align}
&\min \,\left\{
H_{i}(Fx,\,Fy)^{r},\,d_{i}(x,Fx)d_{i}(y,\,Fy)^{r-1},\,d_{i}(y,Fy)^{r}\right\}\nonumber\\[.2pc]
&\qquad\ \ + a_{i}\min \left\{ d_{i}(x,Fy),d_{i}(y,Fx)\right\} \leq [ b_{i}d_{i}(x,Fx)\nonumber\\[.2pc]
&\qquad\ \ + c_{i}d_{i}(x,y)] d_{i}(y,Fy)^{r-1}
\end{align}
for all $i\in I$ and $x,\,y\in X$, where $r\geq 1$ is an integer, $%
a_{i,}\,b_{i},\,c_{i}$ are real numbers such that $0<b_{i}+c_{i}<1$, then $F$
has a fixed point.
\end{theor}

\begin{proof}
Let $x_{0}$ be an arbitrary point in $X$ and consider the sequence
$\{x_{n}\} $ defined by 
\begin{equation*}
x_{1}\in Fx_{0}, x_{2}\in Fx_{1},...,x_{n}\in Fx_{n-1},.... 
\end{equation*}
Let us suppose that $d_{i}(x_{n},\,Fx_{n}) > 0$ for each $i\in I$ and 
$n=0,1,2,...$. (Otherwise for some positive integer $n$, $x_{n}\in
Fx_{n}$ as desired.)

Let $U\in \mathcal{U}$ be an arbitrary entourage. Since $\beta$ is a base
for $\mathcal{U}$, there exists $V(i,\,\varepsilon)\in \beta$ such that 
$V(i,\,\varepsilon)\subseteq U$. Now $y\rightarrow d_{i}(x_{0},\,y)$ is
continuous on the compact set $Fx_{0}$ and this implies that there exists 
$x_{1}\in Fx_{0}$ such that $d_{i}(x_{0},\,x_{1})=d_{i}(x_{0},\,Fx_{0})$.
Similarly, $Fx_{1}$ is compact so there exists $x_{2}\in Fx_{1}$ such that 
$d_{i}(x_{1},\,x_{2}) = d_{i}(x_{1},\,Fx_{1})$. Continuing, we obtain a
sequence $\{x_{n}\}$ such that $x_{n+1}\in Fx_{n}$ and $d_{i}(x_{n},
\,x_{n + 1})=d_{i}(x_{n},\,Fx_{n})$.

For $x=x_{n-1}$, and $y=x_{n}$ by condition (1) we have 
\begin{align*}
&\min \left\{
H_{i}(Fx_{n-1},\,Fx_{n})^{r},\,d_{i}(x_{n-1},\,Fx_{n-1})d_{i}(x_{n},
\,Fx_{n})^{r-1},\,d_{i}(x_{n},\,Fx_{n})^{r}\right\} \\[.2pc]
&\qquad\ \ +a_{i}\,\min \left\{ d_{i}(x_{n-1},\,Fx_{n}),\,d_{i}(x_{n},\,Fx_{n-1})\right\}
 \leq \big[ \,b_{i}\,d_{i}(x_{n-1},Fx_{n-1}) \\[.2pc]
&\qquad\ \ +\,c_{i}\,d_{i}(x_{n-1},
\,x_{n})\big] \,d_{i}(x_{n,\,}Fx_{n})^{r-1}
\end{align*}
or since $d_{i}(x_{n},\,Fx_{n-1})=0,\,x_{n}\in Fx_{n-1}$. Hence we have 
\begin{align*}
&\min \left\{ d_{i}(x_{n},\,x_{n+1})^{r},\,d_{i}(x_{n-1},\,x_{n})d_{i}(x_{n},
\,x_{n+1})^{r-1}\right\} \\[.2pc]
&\quad\ \leq \left[ \,b_{i}\,d_{i}(x_{n-1},x_{n})+\,c_{i}\,d_{i}(x_{n-1},\,x_{n})
\right] \,d_{i}(x_{n,\,}x_{n+1})^{r-1}
\end{align*}
and it follows that 
\begin{align*}
&\min \left\{
d_{i}(x_{n},\,x_{n+1})^{r},\,d_{i}(x_{n-1},\,x_{n})d_{i}(x_{n},
\,x_{n+1})^{r-1}\right\}\\[.2pc]
&\quad\ \leq \,(b_{i}\,+\,c_{i})\,d_{i}(x_{n-1},\,x_{n})\,d_{i}(x_{n,\,}x_{n+1})^{r-1}. 
\end{align*}

Since 
\begin{equation*}
d_{i}(x_{n-1},\,x_{n})\,d_{i}(x_{n},\,x_{n+1})^{r-1}\leq
\,(b_{i}\,+\,c_{i})\,d_{i}(x_{n-1},\,x_{n})\,d_{i}(x_{n,\,}x_{n+1})^{r-1} 
\end{equation*}
is not possible (as $0<b_{i}+c_{i}<1$), we have 
\begin{equation*}
d_{i}(x_{n},\,x_{n+1})^{r}\leq
\,(b_{i}\,+\,c_{i})\,d_{i}(x_{n-1},\,x_{n})\,d_{i}(x_{n,\,}x_{n+1})^{r-1} 
\end{equation*}
or 
\begin{equation*}
\,d_{i}(x_{n},\,x_{n+1})^{r}\leq
\,k_{i}\,\,d_{i}(x_{n-1},\,x_{n})\,d_{i}(x_{n,\,}x_{n+1})^{r-1},
\end{equation*}
where $k_{i} = b_{i} + c_{i\,,\,}0 < k_{i} < 1$.

Proceeding in this manner we get 
\begin{align*}
\,d_{i}(x_{n},\,x_{n+1}) &\leq k_{i}d_{i}(x_{n-1},\,x_{n})\\[.2pc]
&\leq k_{i}^{2}d_{i}(x_{n-2},\,x_{n-1}) \\[.2pc]
&\hskip 1.2pc \vdots \\[.2pc]
&\leq k_{i}^{n}d_{i}(x_{0},\,x_{1}).
\end{align*}

Hence we obtain 
\begin{align*}
d_{i}(x_{n},\,x_{m}) &\leq
d_{i}(x_{n},\,x_{n+1})+\,d_{i}(x_{n+1,\,}x_{n+2})+\cdots+d_{i}(x_{m-1},\,x_{m})\\[.2pc]
&\leq (k_{i}^{n}+k_{i}^{n+1}+\cdots+k_{i}^{m-1})\,\,d_{i}(x_{0},\,x_{1})\\[.2pc]
&\leq k_{i}^{n}(1+k_{i}+\cdots+k_{i}^{m-n-1})\,\,d_{i}(x_{0},\,x_{1}) \\[.2pc]
&\leq \frac{k_{i}^{n}}{1-k}\,\,d_{i}(x_{0},\,x_{1}).
\end{align*}

Since $\underset{n\rightarrow \infty }{\lim }k_{i}^{n}=0$, it follows that
there exists $N(i,\varepsilon)$ such that $d_{i}(x_{n},\,x_{m}) < \varepsilon 
$ and hence $(x_{n},\,x_{m})\in U$ for all $n,m\geq N(i,\varepsilon)$.
Therefore the sequence $\left\{ x_{n}\right\}$ is a Cauchy sequence in the 
$d_{i}$-uniformity on $X$.

Let $S_{p}=\left\{ x_{n}:\,n\geq p\right\} $ for all positive integers $p$
and let $\beta $ be the filter basis $\left\{ S_{p}:p=1,\,2,\,...\right\}$.
Then since $\left\{ x_{n}\right\}$ is a $d_{i}$-Cauchy sequence for each 
$i\in I$, it is easy to see that the filter basis $\beta $ is a Cauchy filter
in the uniform space $\left( X,\,\mathcal{U}\right)$. To see this we first
note that the family $\left\{ V(i,\varepsilon ):i\in I\right\} $ is a base
for $\mathcal{U}$ as $p^{*}=\left\{ d_{i}:i\in I\right\}$. Now since 
$\left\{ x_{n}\right\}$ is a $d_{i}$-Cauchy sequence in $X$, there exists a
positive integer $p$ such that $d_{i}(x_{n},\,x_{m})<\varepsilon$ for 
$m\geq p,\,n\geq p$. This implies that $S_{p}\times S_{p}\subseteq
V(i,\varepsilon)$. Thus given any $U\in \mathcal{U},$ we can find an 
$S_{p}\in \beta $ such that $S_{p}\times S_{p}\subset U$. Hence $\beta $ is a
Cauchy filter in $(X,\,\mathcal{U)}$. Since $\left( X,\mathcal{\,U}
\right)$ is $F$-orbitally\thinspace complete and Hausdorff \thinspace
space, $S_{p}\rightarrow z$ for some $z\in X$. Consequently 
$F(S_{p})\rightarrow Fz$ (follows from the continuity of $F$). Also 
\begin{equation*}
S_{p+1}\subseteq F(S_{p})=\cup \left\{ Fx_{n}:n\geq p\right\} 
\end{equation*}
for $p=1,\,2,\ldots$. It follows that $z\in Fz$. Hence $z$ is a fixed
point of $F$. This completes the proof.\vspace{-.3pc}
\end{proof}

If we take $r=1$ in Theorem $1$, then we obtain the following theorem.

\begin{theor}[\!]
Let $\left( X,\mathcal{\,U}\right)$ be an $F$-orbitally 
complete Hausdorff uniform space defined by
$\left\{ d_{i}:i\in I\right\} = p^{*}$ and $(2^{X},\,\mathcal{U}^{*})$ a
hyperspace{\rm ,} let $F:\,X\rightarrow 2^{X}$ be a continuous mapping and $Fx$
compact for each $x$ in $X$. Assume that 
\begin{align}
&\min \left\{ H_{i}(Fx,\,Fy),\,d_{i}(x,\,Fx),\,d_{i}(y,\,Fy)\right\}\nonumber\\[.2pc]
&\quad\ +a_{i}\min \left\{ \,d_{i}(x,\,Fy),\,d_{i}(y,\,Fx)\right\} \leq
b_{i}\,d_{i}(x,\,Fx)+c_{i}\,d_{i}(x,\,y)
\end{align}
for all $i\in I$ and $x,\,y\in X,$ where $a_{i},\,b_{i},c_{i}$ are real
numbers such that $0 < b_{i} + c_{i} < 1${\rm ,} then $F$ has a fixed point.

We denote that if $F$ is a single valued mapping on $X${\rm ,} then we can write 
$d_{i}(Fx,\,Fy) = H_{i}(Fx,\,Fy), x,\,y\in X,\,i\in I$.
\end{theor}

Thus we obtain the following theorem as a consequence of the Theorem~2.

\begin{theor}[\!]
Let $\left( X,\,\mathcal{U}\right)$ be a $T$-orbitally complete
Hausdorff uniform space and let $T:\,X\rightarrow X$ be a $T$-orbitally
continuous mapping satisfying
\begin{align}
&\min \left\{ d_{i}(Tx,\,Ty),\,d_{i}(x,\,Tx),\,d_{i}(y,\,Ty)\right\}\nonumber\\[.2pc]
&\quad\ +a_{i}\min \left\{ \,d_{i}(x,\,Ty),\,d_{i}(y,\,Tx)\right\} \leq
b_{i}\,d_{i}(x,\,Tx)+c_{i}\,d_{i}(x,\,y)
\end{align}
for all $x,\,y\in X,\,i\in I$ and $a_{i},\,b_{i},c_{i}$ are real numbers
such that $0 < b_{i} + c_{i} < 1$. Then $T$ has a fixed point and which
is unique whenever $a_{i} > c_{i} > 0$.
\end{theor}

\begin{proof}
Define a mapping $F$ of $X$ into $2^{X}$ by putting $Fx=\{Tx\}$ for all 
$x\,$ in $X$. It follows that $F$ satisfies the conditions of Theorem~2.
Hence $T$ has a fixed point.

Now if $a_{i} > c_{i} > 0$, we show that $T$ has a unique fixed point. Assume that 
$T$ has two fixed points $z$ and $w$ which are distinct. Since 
$d_{i}(z,\,Tz) = 0$ and $d_{i}(w,\,Tw) = 0$, then by the condition (2), 
\begin{equation*}
a_{i}\min \left\{ \,d_{i}(z,\,Tw),\,d_{i}(w,\,Tz)\right\} \leq
c_{i}\,d_{i}(z,\,w) 
\end{equation*}
or 
\begin{align*}
a_{i}\,d_{i}(z,\,w) &\leq c_{i}\,d_{i}(z,\,w),\\[.2pc]
d_{i}(z,\,w) &\leq \frac{c_{i}}{a_{i}}\,d_{i}(z,\,w)
\end{align*}
which is impossible. Thus if $a_{i} > c_{i} > 0$, then $T$ has a unique fixed
point in $X$. This completes the proof.\vspace{-.4pc}
\end{proof}

We note that if $a_{i}=-1$ in condition (3), then one gets the following
result as a corollary.

\setcounter{corolla}{3}
\begin{corolla}$\left.\right.$\vspace{.5pc}

\noindent Let $T$ be an orbitally cotinuous self-map of a $T$-orbitally complete
uniform space $\left( X,\mathcal{\,U}\right)$ satisfying the condition
\begin{align*}
&\min \,\left\{ d_{i}(Tx,\,Ty),\,d_{i}(x,\,Tx),\,d_{i}(y,\,Ty)\right\}\\[.2pc]
&\quad\ -\min\,\left\{ d_{i}(x,\,Ty),\,d_{i}(y,\,Tx)\right\} \leq
b_{i}\,d_{i}(x,\,Tx)+c_{i}\,d_{i}(x,\,y),
\end{align*}
$x,\,y\in X,\,i\in I$ and $0 < b_{i} + c_{i} < 1$. Then for each $x\in
X${\rm ,} the sequence $\{T^{n}x\}$ converges to a fixed point of $T$.
\end{corolla}

\begin{rema}
{\rm If we replace the uniform space $(X,\,\mathcal{U})$ in Theorem~3 and
Corollary~4 by a metric space (i.e. a metrizable uniform space), then
Theorem~1 and Corollary~1 of Dhage \cite{2} will follow as special cases of
our results.}
\end{rema}

\section*{Acknowledgement}

This research was supported by the Scientific and Technical Research Council
of Turkey, TBAG-1742 (1999).

\end{document}